\documentclass[11pt]{article}
\usepackage{epsf,amssymb,latexsym,amsmath}
\usepackage{graphicx}
\usepackage{tikz}
\usepackage[autostyle=true]{csquotes}

\newcommand{\proof}{\noindent{\bf Proof.\ }}
\newcommand{\qed}{\hfill $\square$ \bigskip}

\newtheorem{theorem}{\bf Theorem}[section]
\newtheorem{corollary}[theorem]{\bf Corollary}
\newtheorem{lemma}[theorem]{\bf Lemma}
\newtheorem{proposition}[theorem]{\bf Proposition}
\newtheorem{conjecture}[theorem]{\bf Conjecture}

\newtheorem{problem}[theorem]{\bf Problem}
\begin{document}

\title{A counterexample to prism-hamiltonicity of 3-connected planar graphs}

\author{
Simon \v Spacapan\footnote{ University of Maribor, FME, Smetanova 17,
2000 Maribor, Slovenia. e-mail: simon.spacapan @um.si.
}}
\date{\today}

\maketitle

\begin{abstract} 
 The  prism over a graph $G$ is the Cartesian product of $G$ with the complete graph $K_2$. 
A graph $G$ is hamiltonian if there exists a spanning cycle in $G$, and $G$ is prism-hamiltonian  if the prism over $G$ is hamiltonian. 

In [M.~Rosenfeld, D.~Barnette, Hamiltonian circuits in certain prisms, Discrete Math. 5 (1973), 389--394] 
the authors conjectured that every 3-connected planar graph is prism-hamiltonian. 
We construct a counterexample to the conjecture.
\end{abstract}

\noindent
{\bf Key words}: Hamiltonian cycle, planar graph

\bigskip\noindent
{\bf AMS subject classification (2010)}: 05C10, 05C45
% ===================================================================

\section{Introduction}

In 1956 Tutte proved  in \cite{tutte2} that every 4-connected planar graph has a Hamilton cycle. On the other hand 
there exist 3-connected planar graphs that are not hamiltonian. An example of such cubic graph was first found by Tutte in \cite{tutte1}, thereby disproving the 
Tait's conjecture.

A {\em $k$-tree}  is a  tree with maximum degree $k$, and a {\em $k$-walk} in $G$ is a closed walk that visits every vertex of $G$ at most $k$ times. 
It is well known (see for example  \cite{kral} and \cite{bib}), that  the following implications hold. 
\begin{center}
$G$ is hamiltonian $\Rightarrow$ $G$ is traceable $\Rightarrow$     $G$ is prism-hamiltonian      $\Rightarrow$   $G$ has      a spanning 2-walk        $\Rightarrow$ $G$ has a spanning 3-tree
\end{center}
In this hierarchy the property of being hamiltonian is the strongest and existence of a spanning 3-tree is the weakest. 
It was proved recently that  in $P_4$-free graphs prism-hamiltonicity is equivalent to existence of a spanning 2-walk, see \cite{eli}. However for general graphs these properties are 
not equivalent,  moreover they are also not  equivalent in the class of  3-connected planar graphs. 

In \cite{barnette}  it was proved that every 3-connected planar graph has a spanning 3-tree. This was strengthened in \cite{richter}, where 
the authors prove that every 3-connected planar graph has a spanning 2-walk. 

In this paper we address the following conjecture of Rosenfeld and Barnette, see \cite{domneva}, \cite{grunbaum} 
and \cite{kral}\footnote{The conjecture is mentioned in the abstract of \cite{domneva}, and  
in \cite{grunbaum} page 1145.  %In  \cite{domneva} and 
%in \cite{grunbaum} the conjecture is formulated for regular 3-connected planar graphs, although the authors of the conjecture considered it for all 
%3-connected planar  graphs
Rosenfeld mentioned the conjecture during many of his talks \cite{mena}, before it was formulated in \cite{kral} as Conjecture~\ref{domneva}.}.

\begin{conjecture}\label{domneva}
Any 3-connected planar graph is prism-hamiltonian. 
\end{conjecture}

In \cite{domneva} and \cite{grunbaum}  the conjecture is given as a special case of a broader conjecture, claiming that every  graph of a simple 4-polytope is hamiltonian. 
We  construct a counterexample to Conjecture \ref{domneva}.

Many classes of graphs are prism-hamiltonian:  chordal 3-connected  planar graphs 
(also known as {kleetopes}) \cite{kral}, planar near-triangulations  \cite{bib}, Halin graphs \cite{kral}, line graphs of bridgeless graphs \cite{kral},  3-connected cubic graphs \cite{paul} and \cite{kaiser}, and graphs that fulfil special degree conditions \cite{kenta}. A characterization of prism-hamiltonian graphs is also given in \cite{kral2}.
Hamiltonicity of $k$-fold prisms is studied in \cite{domneva}, where the authors prove that for every $k\geq 2$ the $k$-fold prism $G\Box Q_k$
over a 3-connected planar graph $G$
 is hamiltonian (here $Q_k$ denotes the $k$-cube).

We also mention  that some  special subclasses of 3-connected planar graphs are hamiltonian. 
Every triangulation of the plane with at most 3
separating triangles is hamiltonian, see~\cite{jack}. This result was extended recently, in \cite{brink2}, where it is proved that all 3-connected planar graphs with at most three 
3-cuts are hamiltonian.

\section{The counterexample}

%In this section we give a counterexample to Conjecture \ref{domneva}. 
%We start with definiotions and the terminology. 
Let $G=(V(G),E(G))$ and $H=(V(H),E(H))$ be graphs.  The graph $G\cap H$ is the graph with vertex set $V(G)\cap V(H)$ and edge set $E(G)\cap E(H)$,  
and $G\cup H$ is the graph with vertex set $V(G)\cup V(H)$ and edge set $E(G)\cup E(H)$. 
If $S\subseteq V(G)$, then $G-S$ is the graph obtained from $G$ by removing all vertices in $S$, and all edges incident to a vertex in $S$. 
If $S=\{x\}$ is a single vertex, we write $G-x$ instead of $G-\{x\}$. 
 If $P$ is a path in $G$ with endvertices 
$u$ and $v$, then we say that $P$ {\em starts} in $u$ and {\em ends} in $v$, or vice versa. 
A {\em block} of a graph is a maximal connected subgraph without a cutvertex.
Let $K_2$ be the complete graph on two vertices, and let us denote $V(K_2)=\{a,b\}$. For $n\in \mathbb{N}$, we use 
$[n]$ to denote the set of positive integers not larger than $n$.
This notation is used throughout the article.

The {\em Cartesian product} of graphs $G$ and $H$ is the graph with vertex set 
$V(G)\times V(H)$, where vertices $(x_1,y_1)$ and $(x_2,y_2)$ are adjacent in $G\Box H$ if $x_1x_2\in E(G)$ and $y_1=y_2$, or 
$x_1=x_2$ and $y_1y_2\in E(H)$. The {\em prism} over $G$ is the Cartesian product $G\Box K_2$.

A {\em Hamilton path (cycle)} in a graph $G$ is a spanning path  (cycle) in $G$.  
A graph is a {\em hamiltonian} graph if it has a Hamilton cycle, and it is  {\em traceable} if it has a Hamilton path.  

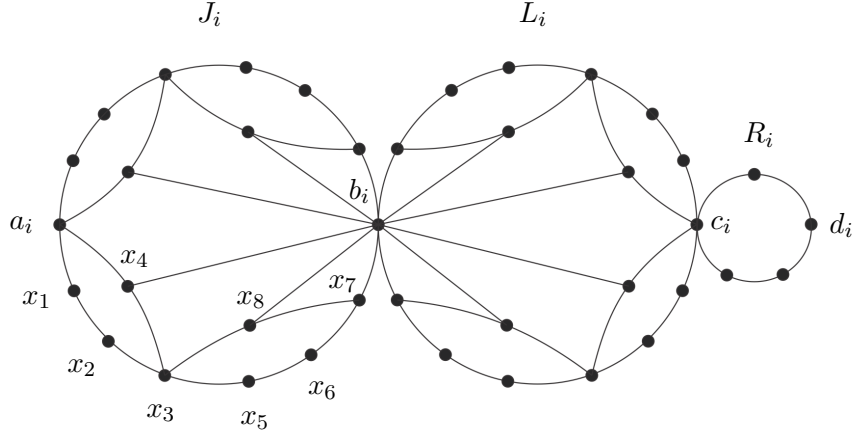
\begin{figure}[bht] 
\begin{center}
\begin{tikzpicture}[scale=1]
\node[inner sep=0pt] at (0,0)
{\includegraphics[width=0.7\textwidth]{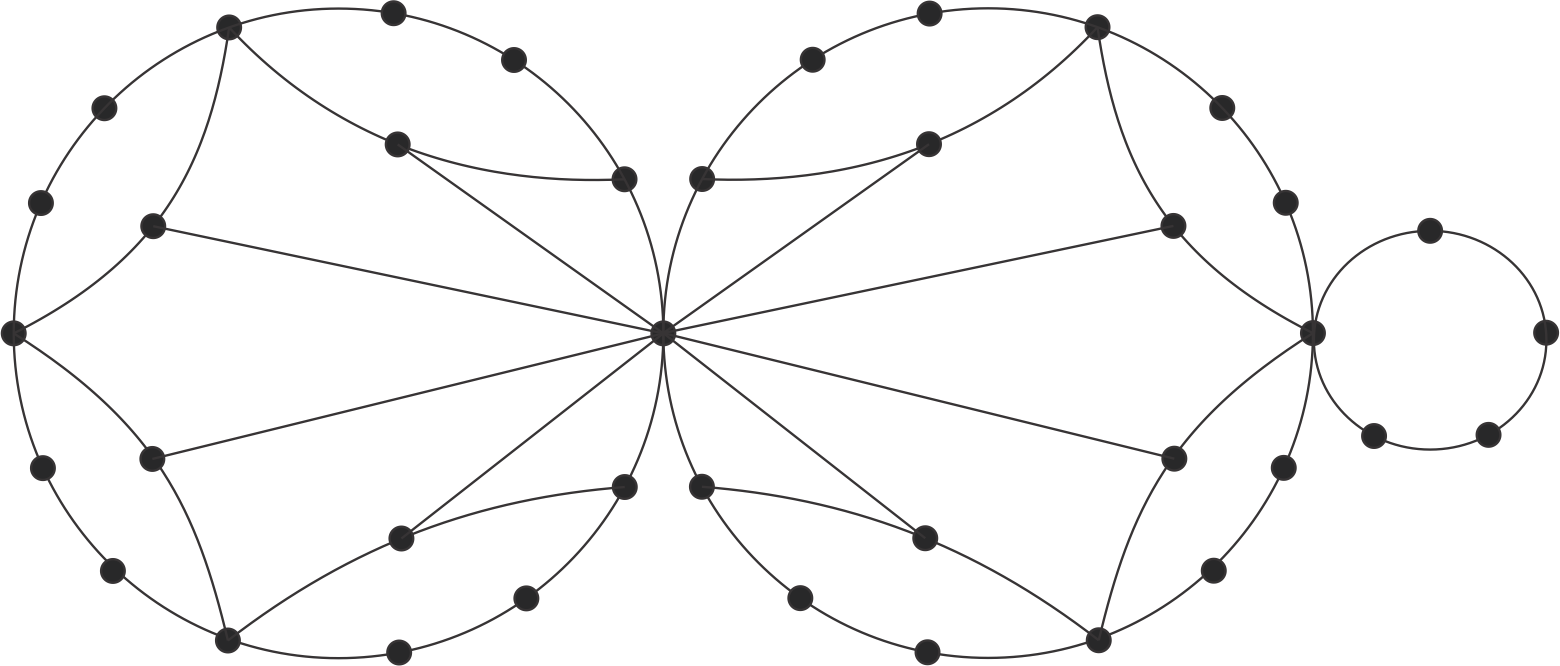}};
%%%%%%%%%%%%%%%%%%%%%%%%%%%%%%%%%%%%%%%zgoraj levo
\path node at (-5.5,0) {$a_i$};
\path node at (-1.01,0.45) {$b_i$};
\path node at (3.8,0) {$c_i$};

\path node at (5.4,0) {$d_i$};

\path node at (-3,2.8) {$J_i$};
\path node at (1.3,2.8) {$L_i$};
\path node at (4.3,1.2) {$R_i$};

%%%%%%%%%%%%%%%%%%%%%%%%%%%%%%%%%%%%%%%%%%%%

\path node at (-5.3,-1) {$x_1$};
\path node at (-4.7,-1.9) {$x_2$};
\path node at (-3.65,-2.5) {$x_3$};
\path node at (-4,-0.45) {$x_4$};
\path node at (-2.4,-2.6) {$x_5$};
\path node at (-2.45,-1) {$x_8$};
\path node at (-1.23,-0.78) {$x_7$};
\path node at (-1.5,-2.2) {$x_6$};

%\path node at (-0.3,-0.8) {$x_9$};
%\path node at (-0.1,-1.98) {$x_{10}$};
%\path node at (0.85,-2.5) {$x_{11}$};
%\path node at (1.2,-1.12) {$x_{13}$};
%\path node at (2.1,-2.5) {$x_{12}$};
%\path node at (2.2,-0.5) {$x_{16}$};
%\path node at (3.62,-1.1) {$x_{15}$};
%\path node at (3.1,-2) {$x_{14}$};

\end{tikzpicture}
\caption{The graph $H_i$.}\label{graf}
\end{center}
\end{figure}
  
Let $H_i$ be the graph shown in Fig.~\ref{graf}. Observe that vertices $b_i$ and $c_i$ are cutvertices of $H_i$, and that $H_i$ has exactly three blocks. 
Let $J_i$ be the block of $H_i$ containing vertex $a_i$, $L_i$ the block containing vertices $b_i$ and $c_i$ and $R_i$ the block of $H_i$ containing 
vertex $d_i$.

Let $Z$ be the subgraph of $J_i$ induced by vertices $a_i,x_1,\ldots,x_8$ (see Fig.~\ref{graf}). Note that $Z$ is the union of two 5-cycles with a common vertex. 
Since $\{(x_3,a),(x_3,b)\}$ is a vertex separator in $Z\Box K_2$, and the cycle containing $a_i$ is odd, we have the following lemma.

\begin{lemma}\label{easy}
 The prism  $Z\Box K_2$ has no Hamilton path with endvertices  $(a_i,a)$ and $(a_i,b)$.
\end{lemma}

\proof Let $C_1$ be the 5-cycle of $Z$ containing vertex $a_i$, and $C_2$ the other 5-cycle.  If $P$ is a path in $Z\Box K_2$ 
 with endvertices  $(a_i,a)$ and $(a_i,b)$ that contains all vertices of $C_2\Box K_2$, then one 
neighbor of $(x_3,a)$ in $P$ is contained in $(C_1-x_3)\Box K_2$ and the other is contained in  $(C_2-x_3)\Box K_2$ (the same is true for vertex $(x_3,b)$). 
The lemma follows form the fact that $C_1$ is odd. \qed

\begin{lemma} \label{osnovna}
 Let $P$ be a path in $J_i\Box K_2$ with the following properties:
\begin{itemize}
\item[(i)] One endvertex of $P$ is $(a_i,a)$ or $(a_i,b)$, and the other is $(b_i,a)$ or $(b_i,b)$ 
\item[(ii)] $P$ contains all vertices of $(J_i- \{a_i,b_i\})\Box K_2$. 
\end{itemize}
Then $(b_i,a)\in P$ and $(b_i,b)\in P$.
\end{lemma}

\proof % Suppose that $P$ is a path that fulfils (i) and (ii).   
Let $P$ be a path as declared in the lemma. Suppose that  $(a_i,a)\notin P$.  Then $P$ starts in $(a_i,b)$. Let $X=\{x_1,\ldots,x_{8}\}$ and $Y=J_i\setminus (X\cup \{a_i,b_i\})$. 
Now there are only two possibilities, either the second vertex of $P$ is contained  in $Y\times V(K_2)$ or it is contained  in $X\times V(K_2)$. Since $J_i$ is symmetric we may assume,  
without loss of generality, that the latter is true (hence the second vertex of $P$ is $(x_{1},b)$ or $(x_{4},b)$).  Since  $(a_i,a)\notin P$, $P$ will eventually go 
through  $(b_i,a)$ or $(b_i,b)$ and enter $Y\times V(K_2)$ to cover vertices of $Y\times V(K_2)$. Since 
$P$ ends with $(b_i,a)$ or $(b_i,b)$ we conclude that $(b_i,a)\in P$ and $(b_i,b)\in P$. 

If  $(a_i,b)\notin P$ the proof is analogous. Assume therefore that $(a_i,a),(a_i,b)\in P$. 
Now, if  $(a_i,a)(a_i,b)\in E(P)$, then the argument is basically the same as above. Either the third vertex of $P$ is contained  in $Y\times V(K_2)$ or it is contained  in $X\times V(K_2)$, with 
the same reasoning as before. So assume that  $(a_i,a)(a_i,b)\notin E(P)$, and assume also that $P$ starts in $(a_i,b)$. Now the only difference of this case is the fact that 
$P$ might enter $X\times V(K_2)$, cover all vertices of $X\times V(K_2)$, and then go through $(a_i,a)$ to enter $Y\times V(K_2)$. 
Since the graph induced by $X\cup\{a_i\}$ is $Z$, this is not possible by Lemma \ref{easy}. \qed

Note that Lemma \ref{osnovna} can be applied to blocks $J_i$ and $L_i$, since they are isomorphic (the isomorphism takes $a_i$ to $c_i$, so these two vertices have the same role when applying Lemma \ref{osnovna}.)

\begin{lemma} \label{osnovna1} 
 There is no path $P$ in $(J_i\cup L_i)\Box K_2$ with the following properties:
\begin{itemize}
\item[(i)] One endvertex of $P$ is  $(a_i,a)$ or $(a_i,b)$, and the other is $(c_i,a)$ or $(c_i,b)$ 
\item[(ii)] $P$ contains all vertices of $((J_i\cup L_i)- \{a_i,c_i\})\Box K_2$. 
\end{itemize}
\end{lemma}

\proof  Suppose on the contrary, that there is such a path $P$.  
Let $P_1=P\cap (J_i\Box K_2)$ and $P_2=P\cap (L_i\Box K_2)$. Then $P_1$ and $P_2$ fulfil the properties (i) and (ii) of Lemma \ref{osnovna}. 
We apply  Lemma \ref{osnovna} to $J_i$ and $L_i$ to find that  $(b_i,a),(b_i,b)\in P_1\cap P_2$. It follows that 
$(b_i,a)(b_i,b)\in E(P)$. We may assume, without loss of generality that 
$(b_i,a)(u,a)\in E(P)$ and $(b_i,b)(v,b)\in E(P)$ for some $u\in J_i$ and $v\in L_i$. 
By symmetry of $J_i$ we may assume (without loss of generality) 
that $u\in Y$ (where $X$ and $Y$ are defined as in Lemma \ref{osnovna}). Since $P_1$ contains all vertices of  $X\times V(K_2)$, we find that $P_1$ will eventually 
go through vertex $(a_i,a)$ or $(a_i,b)$ and enter $X\times V(K_2)$. Then $P_1\cap (Z \Box K_2)$ is a Hamiltonian path in 
$Z \Box K_2$, with endvertices $(a_i,a)$ and $(a_i,b)$. This contradicts Lemma \ref{easy}. 
\qed

It's straightforward to check the following lemma (which is due to the fact that $R_i$ is an odd cycle).

\begin{lemma} \label{osnovna2} 
If $P$ and $P'$ are disjoint paths  in $R_i\Box K_2$  such that 
both paths have one endvertex in  $\{(c_i,a),(c_i,b)\}$ and the other in $\{(d_i,a),(d_i,b)\}$, then
$V(P)\cup V(P')\neq V(R_i\Box K_2)$.
\end{lemma}

Now we define counterexamples to the conjecture. 
Let $P_i$ and $Q_i$ be paths on the boundary walk of $H_i$ from vertex $a_i$ to $d_i$, % (each edge of both paths is incident to the unbounded face of $H_i$), 
and note that there are exactly two such paths, %on the boundary walk of $H_i$, 
the \enquote{upper}  
and the \enquote{lower}. 

%We construct the graph  $G_n$ by gluing graphs $H_i$ together and adding two extra vertices.  
Let $G_n$ be the graph obtained from  graphs $H_i$, where $i\in [n]$, and two additional vertices $x$ and $y$.  To construct the graph $G_n$ we first identify vertex $d_i$ with 
 $a_{i+1}$ for $i\in [n-1]$. 
We may embed the obtained graph in the plane so  that $P_i$ is always the  \enquote{upper}  
path and $Q_i$ the \enquote{lower} path. 
Then we 
draw edges from $x$ to every vertex in  $P_i$ and from $y$ to every vertex in $Q_i$, for $i\in [n]$. Clearly, this can be done so that the obtained graph $G_n$ is a plane graph. 
%It follows that $G_n$ is a planar graph. 

\begin{lemma} \label{planar}   For every $n\in \mathbb{N}$, the graph  $G_n$ is  3-connected.
\end{lemma}

\proof  Suppose that the claim is false and that  $S\subseteq V(G_n)$  is a vertex separator of size 2.  
Note that $\{x,y\}$ is not a vertex separator in $G_n$, therefore we may assume (by symmetry) that $x\notin S$.   

If  $y\in S$, then $S=\{y,u\}$, where $u\neq x$.  Observe that every block of $G_n-\{x,y\}$ contains at most one vertex of $S$ (the vertex $u$). Therefore, 
if $B$ is a block of $G_n-\{x,y\}$, then $B-S$ is connected.   
Since $x$ is adjacent to more than one  vertex of $B$ (in fact $x$ is adjacent to at least three vertices of $B$) 
we find that $(V(B)\cup\{x\})\setminus S$ induces a connected graph.  Since this is true for every block $B$, the graph 
$G_n-S$ is connected. 

If also $y\notin S$, then $S=\{u,v\}$, where $u,v\notin \{x,y\}$. 
If $B-S$ is a connected graph, for every block $B$ of $G_n-\{x,y\}$, the argument is  the same as above. 
Assume therefore that there is a block $B$ of $G_n-\{x,y\}$ such that $B- S$ is disconnected.
This is possible only if $u,v \in V(B)$. It follows that there is at most one block $B$ such that 
$B-S$ is not connected. Note that every connected component of  $B-S$ 
containes an external vertex of $B$ (a vertex incident to the unbounded face of $B$).  It follows that every connected component of $B-S$ is adjacent 
 to $x$ or $y$ in $G_n$. The argument is completed by noting that $x$ and $y$ are adjacent to all other blocks of  $G_n-\{x,y\}$ (by more than one edge), 
and all these blocks remain connected in $G_n-S$.  
\qed

\begin{theorem} If  $n>25$, then $G_n$ is not prism-hamiltonian. %$G_n\Box K_2$ is not Hamiltonian.
\end{theorem}\label{glavni}

%\begin{theorem} If $n$ is large enough, then $G_n$ is not prism-hamiltonian. %$G_n\Box K_2$ is not Hamiltonian.
%\end{theorem}\label{glavni}

\proof
Suppose that  $n>25$, and  that $C$ is a Hamilton cycle in $G_n \Box K_2$.   
There are at most 8 edges of $C$ incident to a vertex in $T=\{(x,a),(x,b)),(y,a),(y,b)\}$. 
%Since every vertex of $G_n-\{x,y\}$ is contained in at most two graphs $H_i$, there are at most 16 graphs $H_i$,
 %such that a vertex  in $H_i\Box K_2$ is adjacent to a vertex in $T$ by an edge of $C$.
Since  $n>25$ there exist consecutive graphs $H_\ell$ and $H_{\ell+1}$, 
such that $C$ has  no edge with one endvertex  in $(H_\ell\cup H_{\ell+1})\Box K_2$ and the other in $T$.

%Note that for every $i\in [n]$, if $u\in \{a_i,b_i,c_i,d_i\}$ and $u\notin \{a_1,d_n\}$, then $\{(u,a),(u,b)\}$ is a vertex separator in $(G_n-\{x,y\})\Box K_2$.  

Let $D_j, j\in [k]$ be connected components of $((H_{\ell}\cup H_{\ell+1})\Box K_2)\cap C$, and let 
$$M=\{(a_{\ell},a),(a_{\ell},b),(d_{\ell+1},a),(d_{\ell+1},b)\}.$$

Each $D_j$ is a path or an isolated vertex. If it is a path, it has both endvertices in $M$. 
If $D_j$ is an isolated vertex, then it is a vertex of $M$ (these properties follow from the assumption that 
 no vertex  in $(H_{\ell}\cup H_{\ell+1})\Box K_2$ is adjacent to a vertex in $T$ by an edge of $C$).  
It follows that $k \leq 3$. Moreover there is either one component $D_j$ which is a path, or there are two such components $D_j$.

Assume the latter, and suppose that $D_1$ and $D_2$ are paths. 
Since $C$ is a Hamilton cycle $D_1$ and $D_2$ contain all vertices of  $(H_{\ell}\cup H_{\ell+1})\Box K_2$.

{\em Case 1.}  Suppose that both, $D_1$ and $D_2$, have one endvertex in $\{ (a_{\ell},a),(a_{\ell},b)\}$ and the other in $\{ (d_{\ell+1},a),(d_{\ell+1},b)\}$.
For every $u\in U=\{b_\ell,c_\ell,d_\ell,b_{\ell+1},c_{\ell+1}\}$ the set  $\{(u,a),(u,b)\}$ is a vertex separator in $(H_{\ell}\cup H_{\ell+1})\Box K_2$.  
It follows that for every $u\in U$ either $(u,a)\in D_1$ and $(u,b)\in D_2$, or $(u,a)\in D_2$ and $(u,b)\in D_1$. 
Therefore $P=D_1\cap (R_{\ell}\Box K_2)$ and   $P'=D_2\cap (R_{\ell}\Box K_2)$ are disjoint paths  that partition  
 the vertex set of the prism $R_{\ell}\Box K_2$. Moreover $P$ and $P'$ have one endvertex in 
 $\{(c_{\ell},a),(c_{\ell},b)\}$ and the other in $\{(d_{\ell},a),(d_{\ell},b)\}$. 
 This contradicts Lemma \ref{osnovna2}.

{\em Case 2.} Suppose that $(a_{\ell},a),(a_{\ell},b)\in D_1$. Then  $(a_{\ell},a),(a_{\ell},b)$ are endvertices of $D_1$,  
and $(d_{\ell+1},a),(d_{\ell+1},b)$ are endvertices of $D_2$. 
If $D_1$ contains a vertex of  $(J_{\ell+1}-\{a_{\ell+1},b_{\ell+1}\})\Box K_2$, 
then $D_1$ contains vertices $(c_{\ell},a), (c_{\ell},b),(d_{\ell},a)$ and $(d_{\ell},b)$.  Therefore $D_1\cap (R_{\ell}\Box K_2)$ is a disjoint union of two paths $P$ and $P'$ that partition  
 the vertex set of the prism $R_{\ell}\Box K_2$. This contradicts Lemma \ref{osnovna2}.
Otherwise $D_2$ contains a vertex of  $(J_{\ell+1}-\{a_{\ell+1},b_{\ell+1}\})\Box K_2$, and therefore $D_2\cap (R_{\ell+1}\Box K_2)$ is a disjoint union of two paths that partition  
 the vertex set of the prism $R_{\ell+1}\Box K_2$, a contradiction. 

Assume now that there is one  component $D_j$ which is a path. Let $D_1$ be the only path, and note that $D_1$ contains all 
vertices of $((H_{\ell}\cup H_{\ell+1})\Box K_2)-M$.  If  $(a_{\ell},a)$ and $(a_{\ell},b)$ are  endvertices of $D_1$, we have 
Case 2 above. 
Otherwise one endvertex of $D_1$ is contained in  $\{ (a_{\ell},a),(a_{\ell},b)\}$ and the other in $\{ (d_{\ell+1},a),(d_{\ell+1},b)\}$.
It follows that $P=D_1\cap((J_\ell\cup L_\ell)\Box K_2)$ is a path that contains all vertices of  $((J_\ell\cup L_\ell)-\{a_\ell,c_\ell\})\Box K_2$. 
Moreover, one endvertex of $P$ is $(a_\ell,a)$ or $(a_\ell,b)$, and the other is $(c_\ell,a)$ or $(c_\ell,b)$. 
This contradicts Lemma \ref{osnovna1}. 
\qed

We note that a slight modification of the above proof  is needed to prove that for sufficiently large $n$,  $G_n\Box K_2$ has no Hamilton path.

\section{Concluding remarks}

 A {\em cactus} is a connected graph $G$ such that every block of $G$  is either a $K_2$ or a cycle. 
An {\em even cactus} is a cactus with no odd cycles. A {\em good cactus} is a cactus such that every vertex is contained in 
at most two blocks.  A {\em  good even cactus} is a cactus that is good and even simultaneously.  A {\em good vertex} of a cactus $G$ is a vertex contained in exactly one block of $G$. 

The authors of \cite{kaiser} considered good even cactuses with the additional property that no vertex is contained in two distinct cycles. 
They prove the following theorem (which is formulated in terms of our terminology). 

\begin{theorem} Every 3-connected cubic graph $G$  has a spanning subgraph $H$, 
such that $H$ is a good even cactus, and every vertex of  $H$ is contained in at most one cycle of $H$.
\end{theorem}

It follows from the above theorem,  and Proposition \ref{bolje}, that 3-connected cubic graphs are 
prism-hamiltonian.
The restriction that no vertex  of a cactus  $G$ is contained in two cycles of $G$ is redundant, when questions of prism-hamiltonicity are addressed. 
This is justified by Proposition \ref{bolje}, and its corollary.

The partitioning result given in \cite{richter} implies that every 3-connected planar graph has a spanning good cactus. 
The existence of a spanning good catctus in $G$ implies the existence of a 2-walk in $G$, while the existence of 
a spanning good even cactus implies prism-hamiltonicity of $G$ - as mentioned in the introduction this is a stronger property. 
The method applied to prove the following proposition is inspired by \cite{bib}. 

\begin{proposition} \label{bolje}
Every good even cactus is prism-hamiltonian.   
\end{proposition}
\proof 
We use induction to prove the following stronger statement. Every prism over a good even cactus $G$ has a Hamilton cycle $C$ such that for  
every good vertex $x$ of $G$, we have $(x,a)(x,b)\in E(C)$.
This is clearly true for all even cycles and $K_2$. Let $G$ be a good even cactus and assume that the statement is true for all good even cactuses with fewer vertices than 
$|V(G)|$. If all vertices of $G$ are good, then $G$ is an even cycle or $K_2$. Otherwise, there  is a  vertex $u$, which is not a good vertex of $G$. 
Hence, $u$ is contained in exactly two blocks of $G$.

Let $G_1'$ and $G_2'$ be connected components of $G-x$, and let $G_1=G-G_2'$ and $G_2=G-G_1'$.
Both, $G_1$ and $G_2$, are good even cactuses. Moreover, $x$ is a good vertex in $G_i$, for $i=1,2$. 
By  induction hypothesis there is a Hamilton cycle $C_i$ n $G_i$ such that $C_i$ uses the edge $e=(x,a)(x,b)$ in $G_i$. 
The desired Hamilton cycle in $G$ is $(C_1\cup C_2)-e$. Observe that every good vertex of $G$ is a good vertex of $G_1$ or $G_2$. 
It follows that for  
every good vertex $x$ of $G$, we have $(x,a)(x,b)\in E(C)$.
\qed

\begin{corollary}
Every graph $G$, that has  a good even cactus  $H$ as a spanning subgraph,  is prism-hamiltonian. 
\end{corollary}

The counterexample to prism-hamiltonicity of 3-connected planar graphs, given in Theorem \ref{glavni}, was constructed via a construction of a
 graph with no spanning good even cactus. We conclude this article with open problems. 

%the question, if 
%\enquote{a 3-connected planar graph $G$  is prism-hamiltonian if and only if $G$ has a good even cactus $H$ as a spanning subgraph} ? 

\begin{problem}
Prove or disprove the following statement. If a 3-connected planar graph $G$  is prism-hamiltonian, then $G$ has a good even cactus $H$ as a spanning subgraph. 
\end{problem}

The following problem is due to Rosenfeld, in fact it's conjectured that question (1) has a positive answer \cite{mena}. 

\begin{problem}
(1) Is every 4-connected  4-regular graph prism-hamiltonian ? (2) Is also every 3-connected  4-regular graph prism-hamiltonian ?  
\end{problem}

%\begin{problem}
%Is every 3-connected  regular planar graph prism-hamiltonian ? 
%\end{problem}

\noindent {\bf Acknowledgement:} 
The author thanks M.~Rosenfeld for helpful  comments on the origin of Conjecture~\ref{domneva}. This work was supported by the Ministry of Education of Slovenia  [grant numbers P1-0297, J1-9109].

\end{document}